\documentclass{amsproc}

\usepackage[colorlinks=true,hyperindex, linkcolor=magenta, pagebackref=false, citecolor=cyan,pdfpagelabels]{hyperref}
\hypersetup{linktocpage}
\usepackage{amssymb}
\usepackage{amsmath}
\usepackage[matrix,arrow,curve,cmtip]{xy}
\entrymodifiers={+!!<0pt,\fontdimen22\textfont2>} 
\usepackage{units}
\usepackage{enumitem}
\setenumerate[1]{leftmargin=2em}
\usepackage{eucal}
\usepackage[nameinlink,capitalise,noabbrev]{cleveref} 

\numberwithin{equation}{section}

\theoremstyle{plain}   

\theoremstyle{definition}

\theoremstyle{remark}

\DeclareMathOperator{\id}{id}

\DeclareMathOperator{\can}{can}
\DeclareMathOperator{\Spec}{Spec}
\DeclareMathOperator{\Proj}{Proj}
\DeclareMathOperator{\Aut}{Aut}
\DeclareMathOperator{\QCoh}{QCoh}
\DeclareMathOperator{\Fil}{Fil}
\DeclareMathOperator{\gr}{gr}
\DeclareMathOperator{\dR}{dR}
\DeclareMathOperator{\Sp}{Sp}

\DeclareMathOperator{\HH}{HH}
\DeclareMathOperator{\HC}{HC}
\DeclareMathOperator{\HP}{HP}
\DeclareMathOperator{\THH}{THH}
\DeclareMathOperator{\TC}{TC}
\DeclareMathOperator{\TP}{TP}
\DeclareMathOperator{\CAlg}{CAlg}
\DeclareMathOperator{\CycSp}{CycSp}

\begin{document}

\title{Topological cyclic homology and the Fargues--Fontaine curve}

\author{Lars Hesselholt}
\address{Nagoya University, Japan, and University of Copenhagen, Denmark}
\email{larsh@math.nagoya-u.ac.jp}

\thanks{The author was partially supported by the Danish National
  Research Foundation through the Copenhagen Center for Geometry and
  Topology (DNRF151) and by JSPS Grant-in-Aid for Scientific Research
  number 21K03161.}

\maketitle

\section*{Introduction}

This paper is an elaboration of my lecture at the conference. The
purpose is to explain how the Fargues--Fontaine curve and its
decomposition into a punctured curve and the formal neighborhood of
the puncture naturally appear from various forms of topological
cyclic homology and maps between them. I make no claim of
originality. My purpose here is to highlight some of the spectacular
material contained in the papers of Nikolaus--Scholze~\cite{nikolausscholze},
Bhatt--Morrow--Scholze~\cite{bhattmorrowscholze2}, and
Antieau--Mathew--Morrow--Nikolaus~\cite{antieaumathewmorrownikolaus}
on topological cyclic homology and in the book by
Fargues--Fontaine~\cite{farguesfontaine} on their revolutionary
curve.

\section{The Fargues--Fontaine curve}

We give a brief introduction to the Fargues--Fontaine curve and refer
to their book~\cite{farguesfontaine} for details. The lecture
notes by Lurie~\cite{lurieff} are also helpful.

We define a completely valued field to be a pair $(C,\mathcal{O}_C)$
of a field $C$ and a proper subring $\mathcal{O}_C \subset C$ that is
a complete valuation ring of rank $1$. The field
$C$ is the quotient field of $\mathcal{O}_C$ and the requirement that
the inclusion $\mathcal{O}_C \subset C$ be proper is equivalent to the
requirement that the value group of $\mathcal{O}_C$ be non-trivial. An
isomorphism of completely valued fields $\psi \colon (C,\mathcal{O}_C) \to
(C',\mathcal{O}_{C'})$ is an isomorphism of fields $\psi \colon C \to
C'$ that restricts to an isomorphism of rings $\psi \colon
\mathcal{O}_C \to \mathcal{O}_{C'}$.

We will be interested in the case, where $C$ is algebraically closed
and $\mathcal{O}_C$ is of residue characteristic $p > 0$. 
In this case, ``tilting'' is a correspondence that to an algebraically
closed completely valued field $(C,\mathcal{O}_C)$ of residue
characteristic $p > 0$ assigns the algebraically closed completely valued
field $(C^{\flat},\mathcal{O}_{C^{\flat}})$ of characteristic $p > 0$,
where $\mathcal{O}_{C^{\flat}}$ is the limit of the diagram of
$\mathbb{F}_p$-algebras
$$\xymatrix@C=7mm{
  { \cdots } \ar[r]^-{\varphi} &
  { \mathcal{O}_C/p } \ar[r]^-{\varphi} &
  { \mathcal{O}_C/p } \ar[r]^-{\varphi} &
  { \mathcal{O}_C/p } \cr
}$$
with $\varphi$ the Frobenius. The non-trivial fact that $C^{\flat}$ is
algebraically closed is proved
in~\cite[Proposition~2.1.11]{farguesfontaine}. 
The tilting correspondence is many-to-one. More
precisely, given an algebraically closed completely valued field
$(F,\mathcal{O}_F)$ of characteristic $p > 0$, an untilt of
$(F,\mathcal{O}_F)$ is a pair
$((C,\mathcal{O}_C),\iota)$ of an algebraically closed completely
valued field $(C,\mathcal{O}_C)$ and an isomorphism of completely
valued fields
$$\xymatrix{
  { (F,\mathcal{O}_F) } \ar[r]^-{\iota} &
  { (C^{\flat},\mathcal{O}_{C^{\flat}}). } \cr
}$$
The isomorphism $\iota$ gives rise to an isomorphism of value groups
$$\xymatrix{
{ F^{\times}/\mathcal{O}_F^{\times} } \ar[r] &
{ C^{\times}/\mathcal{O}_C^{\times} } \cr
}$$
and the ``perfectoid'' nature of the respective valuation rings
implies that for every pair of elements $\varpi_F \in \mathcal{O}_F$
and $\varpi_C \in \mathcal{O}_C$ of equal and sufficiently small but
nonzero valuation, the isomorphism $\iota$ gives rise to an
isomorphism of rings
$$\xymatrix{
{ \mathcal{O}_F/\varpi_F } \ar[r] &
{ \mathcal{O}_C/\varpi_C. } \cr
}$$
We note that the common ring is not a field, but rather is a highly
non-noetherian nilpotent thickening of its (algebraically closed)
residue field. 

An isomorphism of untilts $\psi
\colon ((C,\mathcal{O}_C),\iota) \to ((C',\mathcal{O}_{C'},\iota')$ is
an isomorphism of completely valued fields $\psi \colon
(C,\mathcal{O}_C) \to (C',\mathcal{O}_{C'})$ with the property that
the diagram
$$\begin{xy}
  (0,0)*+{ (F,\mathcal{O}_F) }="1";
  (25,7)*+{ (C^{\flat},\mathcal{O}_{C^{\flat}}) }="2";
  (25,-7)*+{ (C'{}^{\flat},\mathcal{O}_{C'{}^{\flat}}) }="3";
  { \ar@/^.5pc/^-(.55){\iota} "2";"1";};
  { \ar@/_.5pc/^-(.52){\iota'} "3";"1";};
  { \ar^-{\psi^{\flat}} "3";"2";};
\end{xy}$$
commutes. Up to isomorphism, there is one untilt of
$(F,\mathcal{O}_F)$ of characteristic $p$. But there are many
non-isomorphic untilts of $(F,\mathcal{O}_F)$ of characteristic $0$.
We note that if $((C,\mathcal{O}_C),\iota)$ is an untilt of
$(F,\mathcal{O}_F)$, then so is $((C,\mathcal{O}_C),\iota \circ
\varphi)$, and that if $C$ is of characteristic $0$, then these two
untilts are non-isomorphic. So the group $\varphi^{\mathbb{Z}}$ acts
freely on the set of isomorphism classes of characteristic $0$ untilts
of $F$. Fargues and Fontaine show
in~\cite[Th\'{e}or\`{e}me~6.5.2]{farguesfontaine} that their curve
$$\xymatrix{
  { X = X_F } \ar[r] &
  { \Spec(\mathbb{Q}_p) } \cr
}$$
parametrizes the orbits of this action in the sense that there is a canonical bijection from the set $|X|$ of
closed points onto the set of
$\varphi^{\mathbb{Z}}$-orbits of isomorphism classes of characteristic
$0$ untilts of
$(F,\mathcal{O}_F)$. Moreover, if $x \in |X|$ corresponds to the class
of $((C,\mathcal{O}_C),\iota)$ via this bijection, then the residue field $k(x)$ is
canonically isomorphic to $C$. The additional structure of the
valuation ring $\mathcal{O}_C \subset C$ and the isomorphism
$\iota \colon (F,\mathcal{O}_F) \to (C^{\flat},\mathcal{O}_{C^{\flat}})$ arises from the
``perfectoid'' nature of the situation, a fact that is hidden away in
the proof of~\cite[Th\'{e}or\`{e}me~6.5.2]{farguesfontaine}.

One may wonder if all untilts of $(F,\mathcal{O}_F)$ have abstractly
isomorphic underlying completely valued fields. More precisely, given
an untilt $((C,\mathcal{O}_C),\iota)$, we let
$$\xymatrix{
  { \Aut(F,\mathcal{O}_F)/(\Aut(C,\mathcal{O}_C) \times
    \varphi^{\mathbb{Z}}) } \ar[r] &
  { |X| } \cr
}$$
be the map that to $\psi$ assigns the class of
$((C,\mathcal{O}_C),\iota \circ \psi)$ and ask if this map is a
bijection. It follows from~\cite[Corollaire~2.2.23]{farguesfontaine}
and~\cite[Corollary~6]{poonen} that the answer to this question is
affirmative if and only if 
the valued field $(F,\mathcal{O}_F)$ is spherically complete (a.k.a.\
maximally complete) in the sense that every decreasing sequence of
discs in $F$ has non-empty intersection. This is a stronger condition than
completeness, which is the condition that every decreasing sequence of
discs, whose radii tend to $0$, has non-empty intersection. For
example, the completion $\mathbb{C}_p$ of an algebraic closure of
$\mathbb{Q}_p$ with respect to the unique extention of the $p$-adic
absolute value is not spherically complete, and neither is its tilt
$\mathbb{C}_p^{\flat}$.

The Fargues--Fontaine curve behaves as a proper, regular,
connected curve of genus $0$ over $\mathbb{Q}_p$, except that it is
not of finite type. Strikingly, the map 
$$\xymatrix{
  { \mathbb{Q}_p } \ar[r] &
  { H^0(X,\mathcal{O}_X) } \cr
}$$
induced by the structure map is an isomorphism. Hence, there is a
great discrepancy between the field $\mathbb{Q}_p$ of global sections
and the residue fields $k(x)$ at closed points $x \in |X|$, all of
which all are algebraically closed completely valued extensions of
$\mathbb{Q}_p$ by~\cite[Th\'{e}or\`{e}me~6.5.2]{farguesfontaine}. This
phenomenon is one that the Fargues--Fontaine curve shares with the
twistor projective line.  The complex projective line
$\mathbb{P}_{\mathbb{C}}^1$ has two real forms, namely, the real
projective line $\mathbb{P}_{\mathbb{R}}^1$ and the twistor projective
line
$$\xymatrix{
  { X = \widetilde{\mathbb{P}}_{\mathbb{R}}^1 } \ar[r] &
  { \Spec(\mathbb{R}), } \cr
}$$
which is the Brauer--Severi variety of the quaternions
$\mathbb{H}$. Its ring of global sections is $\mathbb{R}$, whereas the
residue field at every closed point $x \in |X|$ is an algebraically
closed field that contains $\mathbb{R}$. These residue fields are of
course all isomorphic to $\mathbb{C}$. In particular, the twistor
projective line has no real points.

\section{Curves}

In general, a connected curve over a field 
$$\xymatrix{
  { X } \ar[r] &
  { \Spec(k) } \cr
}$$
can be understood as follows. If we choose a regular closed point
$\infty \in |X|$, then the open complement $X \smallsetminus
\{\infty\} \subset X$ is affine, and the quasi-coherent ideal
$$\xymatrix{
  { 0 } \ar[r] &
  { \mathcal{O}_X(-1) } \ar[r] &
  { \mathcal{O}_X } \ar[r] &
  { i_{\infty*}k(\infty) } \ar[r] &
  { 0 } \cr
}$$
is invertible. Moreover, writing $\mathcal{O}_X(n)$ for its
$(-n)$-fold tensor power, we have
$$\xymatrix{
  { X \simeq \Proj(P) } \ar[r] &
  { \Spec(k), } \cr
}$$
where $P$ is the graded $k$-algebra
$$\xymatrix{ P = \bigoplus_{n \in \mathbb{Z}}
  H^0(X,\mathcal{O}_X(n)). }$$
Associated with the closed point $\infty \in |X|$, we have the maps
$$\xymatrix@C=10mm{
  { U } \ar[r]^-{\phantom{j}i\phantom{j}} &
  { X } &
  { Y } \ar[l]_-{j} \cr
}$$
of formal schemes over $k$, where $i$ and $j$ are the
canonical maps from the open complement of $\{\infty\} \subset X$ and
the formal completion along $\{\infty\} \subset X$, respectively. In
Grothendieck's philosophy, the map $j$ should be viewed as the 
open and affine immersion of a tubular neighborhood of $\{\infty\} \subset
X$, whereas the map $i$ should be viewed as the closed immersion of
the closed complement of said neighborhood. This point of view is
substantiated by the stable recollement
\vspace{-1mm}
$$\begin{xy}
  (0,0)*+{ \QCoh(U) };
  (30.5,0)*+{ \QCoh(X) };
  (61,0)*+{ \QCoh(Y) };
  (7,0)*+{}="1";
  (23,0)*+{}="2";
  (38,0)*+{}="3";
  (54,0)*+{}="4";
  (0,-5)*+{};
  { \ar@<.5ex>@/^1pc/_-{i^!} "1";"2";};
  { \ar^-{i_* \simeq\, i_!} "2";"1";};
  { \ar@<-.9ex>@/_1pc/_-{i^*} "1";"2";};  
  { \ar@<.5ex>@/^1pc/_-{j_*} "3";"4";};
  { \ar^-{j^! \simeq\, j^*} "4";"3";};
  { \ar@<-.9ex>@/_1pc/_-{j_!} "3";"4";};  
\end{xy}$$
among the corresponding stable $\infty$-categories of quasi-coherent
modules.\footnote{\,Let $\smash{ Y^{(n)} \subset X }$ be the $n$th
  infinitesimal neighborhood of $\{\infty\} \subset X$. We define
  $\QCoh(Y)$ to be the colimit $\smash{
    \operatorname{colim}_n\QCoh(Y^{(n)}) }$ in the $\infty$-category
  of presentable $\infty$-categories and right adjoint functors, or
  equivalently, as the limit $\smash{ \lim_n\QCoh(Y^{(n)}) }$ in the
  $\infty$-category of presentable $\infty$-categories and left
  adjoint functors. The functor $\smash{ j^! \simeq j^*}$ is Grothendieck's
  local cohomology.} We give a proof of this is
in~\cite{hesselholtpstragowski}, but we also refer
to~\cite[Lecture~V]{clausenscholzecomplex} for an explanation of the
open-closed reversal in this situation. So we have a cartesian diagram
in the $\infty$-category $\QCoh(X)$ of quasi-coherent
$\mathcal{O}_X$-modules
$$\xymatrix{
  { \mathcal{O}_X(n) } \ar[r] \ar[d] &
  { i_*i^*\mathcal{O}_X(n) } \ar[d] \cr
  { j_*j^*\mathcal{O}_X(n) } \ar[r] &
  { i_*i^*j_*j^*\mathcal{O}_X(n). } \cr
}$$
We interpret $\mathcal{O}_X(n)$ as the sheaf of meromorphic functions
on $X$ that are regular away from $\infty$ and whose pole order at
$\infty$ is at most $n$, and $i_*i^*\mathcal{O}_X(n)$ as the sheaf of 
meromorphic functions on $X$ that are regular away from
$\infty$. Similarly, we interpret $j_*j^*\mathcal{O}_X(n)$ as the sheaf of 
formal meromorphic functions near $\infty$ that are regular away from
$\infty$ and whose pole order at $\infty$ is at most $n$, and
$i_*i^*j_*j^*\mathcal{O}_X(n)$ as the sheaf of formal meromorphic
functions near $\infty$ that are regular away from $\infty$. It
follows that the $k$-vector spaces $H^0(X,\mathcal{O}_X(n))$ and
$H^1(X,\mathcal{O}_X(n))$ are canonically identified with the limit
and colimit, respectively, of the diagram
\begin{equation}\label{eq:diagram}
  \xymatrix{
    {} &
    { H^0(X, i_*i^*\mathcal{O}_X(n)) } \ar[d] \cr
    { H^0(X,j_*j^*\mathcal{O}_X(n)) } \ar[r] &
    { H^0(X,i_*i^*j_*j^*\mathcal{O}_X(n)). } \cr
  }
\end{equation}
As a simple example, let us first consider the complex projective line
$$\xymatrix{
  { X = \mathbb{P}_{\mathbb{C}}^1 } \ar[r] &
  { \Spec(\mathbb{C}). } \cr
}$$
If we choose a coordinate $t^{-1}$ at $\infty \in |X|$, then the
diagram~(\ref{eq:diagram}) takes the form
\begin{equation}\label{eq:projectiveline}
  \xymatrix@C=10mm{
    {} &
    { \mathbb{C}[t] } \ar[d] \cr
    { t^n\mathbb{C}[[t^{-1}]] } \ar[r] &
    { \mathbb{C}((t^{-1})) } \cr
  }
\end{equation}
with the two maps given by the canonical inclusions. So for $n \geq
0$, we have
$$\begin{aligned}
  H^0(X,\mathcal{O}_X(n)) & = \mathbb{C} \cdot \{1,t,\dots,t^n\} \cr
  H^1(X,\mathcal{O}_X(n)) & = 0, \cr
\end{aligned}$$
where $\mathbb{C} \cdot S$ is the $\mathbb{C}$-vector space generated by
$S$, whereas for $n < 0$, we have
$$\begin{aligned}
  H^0(X,\mathcal{O}_X(n)) & = 0 \cr
  H^1(X,\mathcal{O}_X(n)) & = \mathbb{C} \cdot \{t^{-1},\dots,t^{n+1}\}. \cr
\end{aligned}$$
In particular, we find that the graded $\mathbb{C}$-algebra $P$ is
given by
$$\xymatrix{
  { P \simeq \bigoplus_{n \geq 0} \mathbb{C} \cdot \{1,t,\dots,t^n\} } \ar[r] &
  { \mathbb{C}[x,y], } \cr
}$$
where the $n$th component of the right-hand isomorphism takes $t^i$ to $x^iy^{n-i}$.

In the case of the twistor projective line
$$\xymatrix{
  { X = \widetilde{\mathbb{P}}_{\mathbb{R}}^1 } \ar[r] &
  { \Spec(\mathbb{R}), } \cr
}$$
the diagram~(\ref{eq:diagram}) takes the form
\begin{equation}\label{eq:twistorprojectiveline}
  \xymatrix{
    {} &
    { \mathbb{R}[u,v]/(u^2+v^2+1) } \ar[d] \cr
    { t^n\mathbb{C}[[t^{-1}]] } \ar[r] &
    { \mathbb{C}((t^{-1})) } \cr
  }
\end{equation}    
with the horizontal map given by the canonical inclusion and with the
vertical map obtained by solving the complex linear system of equations
$$\begin{aligned}
  u + iv & = t \cr
  u - iv & = -t^{-1}. \cr
\end{aligned}$$
We remark that now the identification of the terms in the bottom row
uses the fact that a complete discrete valuation ring of residue
characteristic $0$ admits a unique coefficient field, which, in
general, relies on the axiom of choice. Thus, the $\mathbb{C}$-algebra 
structure on the two rings is somewhat of a mirage, as opposed to the
$\mathbb{C}$-algebra structure on the associated graded rings for the
$t$-adic filtrations. This
filtration, in turn, induces a $\mathbb{Z}$-graded ascending
filtration
$$\cdots \subset \Fil_{m-1}H^i(X,\mathcal{O}_X(n)) \subset
\Fil_mH^i(X,\mathcal{O}_X(n)) \subset \cdots$$
of the cohomology groups in question, and we find that for $n \geq 0$,
$$\gr_mH^0(X,\mathcal{O}_X(n)) \simeq \begin{cases}
  \mathbb{R} \cdot 1 & \text{if $m = 0$} \cr
  \mathbb{C} \cdot t^m & \text{if $0 < m \leq n$} \cr
\end{cases}$$
are the only nonzero graded pieces, whereas for $n < 0$,
$$\gr_mH^1(X,\mathcal{O}_X(n)) \simeq \begin{cases}
  \mathbb{C} \cdot t^m & \text{if $n+1 \leq m \leq -1$} \cr
  \mathbb{C} / \hspace{1pt}\mathbb{R} \cdot 1 & \text{if $m = 0$} \cr
\end{cases}$$
are the only nonzero graded pieces. In particular, the
quotient $\mathbb{C}/\hspace{1pt}\mathbb{R}$ of the residue field at
$\infty$ by the subfield of global sections appears as
$H^1(X,\mathcal{O}_X(-1))$. We find that there is an isomorphism of
graded $\mathbb{R}$-algebras
$$\xymatrix{
  { P = \bigoplus_{n \geq 0}H^0(X,\mathcal{O}_X(n)) } \ar[r] &
  { \mathbb{R}[x,y,z]/(x^2+y^2+z^2) } \cr
}$$
whose $n$th component takes the class of $u^iv^j$ to
$x^iy^jz^{n-(i+j)}$. 

In the case of the Fargues--Fontaine curve
$$\xymatrix{
  { X = X_F } \ar[r] &
  { \Spec(\mathbb{Q}_p), } \cr
}$$
the diagram~(\ref{eq:diagram}) is expressed in terms of Fontaine's period
rings\footnote{\,The notation $B_e$ stems
  from~\cite[(3.7.2)]{blochkatotamagawa}. The subscript $e$ indicates
  the exponential case.} as
\begin{equation}\label{eq:farguesfontainecurve}
  \xymatrix@C=10mm{
    {} &
    { B_e } \ar[d] \cr
    { \Fil_nB_{\dR} } \ar[r] &
    { B_{\dR} } \cr
  }
\end{equation}
The ring $B_{\dR}$ is a complete discrete valuation field with residue
field $C = k(\infty)$,\footnote{\,While there does exist an isomorphism
$B_{\dR} \simeq C((t^{-1}))$ of complete discrete valuation fields,
such an isomorphism is highly non-canonical and not useful in
practice.} and the horizontal map is the canonical
inclusion of the $(-n)$th power of the maximal ideal of its valuation
ring $B_{\dR}^+ = \Fil_0B_{\dR}$. The ring $B_e$ is a principal ideal
domain, the nonzero ideals of which are in canonical one-to-one
correspondence with the closed points $x \in |X \smallsetminus
\{\infty\}|$. This non-trivial fact, instigated by Berger's
discovery that every finitely generated ideal in $B_e$ is
principal~\cite[Proposition~1.1.9]{berger}, was the discovery which led Fargues and 
Fontaine to realize that they had a curve in their hands. The proof is
given in~\cite[Th\'{e}or\`{e}me~6.5.2]{farguesfontaine}, and Colmez'
recounting of this discovery process
in~\cite[Pr\'{e}face]{farguesfontaine} is quite illuminating. The
discrete valuation on $B_{\dR}$ again gives rise to an ascending
$\mathbb{Z}$-graded filtration
$$\cdots \subset \Fil_{m-1}H^i(X,\mathcal{O}_X(n)) \subset
\Fil_mH^i(X,\mathcal{O}_X(n)) \subset \cdots$$
of the cohomology groups, and for $n \geq 0$,
$$\gr_mH^0(X,\mathcal{O}_X(n)) = \begin{cases}
  \mathbb{Q}_p \cdot 1 & \text{if $m = 0$} \cr
  C \cdot t^m & \text{if $0 < m \leq n$} \cr
\end{cases}$$
are the only nonzero graded pieces, whereas for $n < 0$,
$$\gr_mH^1(X,\mathcal{O}_X(n)) = \begin{cases}
  C \cdot t^m & \text{if $n+1 \leq m \leq -1$} \cr
  C / \hspace{.8pt}\mathbb{Q}_p \cdot 1 & \text{if $m = 0$} \cr
\end{cases}$$
are the only nonzero graded pieces. Here $t^{-1} \in \Fil_{-1}B_{\dR}$
is a local parameter at $\infty$, that is, a generator of this
$B_{\dR}^+$-module, which is free of rank $1$. Again, the quotient
$C/\hspace{.8pt}\mathbb{Q}_p$ of the residue field at $\infty$ and the
subfield of global sections appears as $H^1(X,\mathcal{O}_X(-1))$, but
this now an infinite dimensional $\mathbb{Q}_p$-vector space,
reflecting the fact that the Fargues--Fontaine curve is not of finite
type over $\mathbb{Q}_p$. Nevertheless, Fargues and Fontaine show
in~\cite[Th\'{e}or\`{e}me~6.2.1]{farguesfontaine} that the graded
$\mathbb{Q}_p$-algebra
$$\textstyle{ P = \bigoplus_{n \geq 0} H^0(X,\mathcal{O}_X(n)) }$$
is a graded unique factorization domain, all of whose irreducible
elements are of degree $1$. Thus, closed points $x \in |X|$ are in
canonical one-to-one correspondence with $\mathbb{Q}_p$-lines in
$H^0(X,\mathcal{O}_X(1))$, or equivalently, extensions of the form
$$\xymatrix{
  { 0 } \ar[r] &
  { \mathbb{Q}_p } \ar[r] &
  { H^0(X,\mathcal{O}_X(1)) } \ar[r] &
  { k(x) } \ar[r] &
  { 0; } \cr
}$$
see~\cite[Th\'{e}or\`{e}me~6.5.2]{farguesfontaine}.

\section{The formal neighborhood of $\infty \in |X|$}

We proceed to explain how to obtain the
diagram~(\ref{eq:farguesfontainecurve}) for the Fargues--Fontaine
curve from a diagram of spectra. We begin with the lower line in the
diagram, and recall some general theory
following Nikolaus--Scholze~\cite[Theorem~I.4.1]{nikolausscholze}.

We will use the language of $\infty$-categories following Joyal and
Lurie~\cite{luriehtt}, but we will use the term ``anima'' or
``animated set'' for what Lurie calls a space to emphasize that these
should be viewed rather as sets with internal symmetries rather than
anything resembling a topological space.

If $f \colon T \to S$ is any map of anima, then we have the
restriction along $f$, $f^! \simeq f^*$, and the left and right Kan
extensions along $f$, $f_!$ and $f_*$, between the corresponding
$\infty$-categories of spectrum-valued presheaves:
\vspace{-1mm}
$$\begin{xy}
  (0,0)*+{ \Sp^T };
  (23,0)*+{ \Sp^S };
  (3,0)*+{}="1";
  (20,0)*+{}="2";
  (0,-5)*+{};
  { \ar@<.9ex>@/^.9pc/^-{f_!} "2";"1";};
  { \ar_-{f^! \simeq\, f^*} "1";"2";};
  { \ar@<-.5ex>@/_.9pc/^-{f_*} "2";"1";};
\end{xy}$$
In fact, these functors are part of a six-functor
formalism on the $\infty$-category of anima in the sense of
Mann~\cite[Definition~A.5.7]{mann}. Indeed, this claim is a trivial
instance of~\cite[Proposition~A.5.10]{mann}, where every map of anima
$f \colon T \to S$ is declared to be a local isomorphism, and where a
map of anima $f \colon T \to S$ is declared to be proper if its fibers
are compact projective anima, or equivalently, anima that are equivalent
to finite sets, including the empty set. In this six-functor
formalism, we have $f^! \simeq f^*$ for every map
$f \colon T \to S$, which reflects the fact that every map of anima is
a local isomorphism. Now, the $\infty$-categories $\Sp^S$ and $\Sp^T$
are compactly 
generated  presentable stable $\infty$-categories, and the
``homology'' functor $f_!$ preserves compact objects, because its
right adjoint $f^!$ preserves all colimits. However, the
``cohomology'' functor $f_*$ does not preserve compact objects, and
by~\cite[Theorem~I.4.1]{nikolausscholze}, there is a unique a map $f_*
\to f_*^T$ to a ``Tate cohomology'' functor that takes all compact
objects to zero and that is initial with this property. The fiber of
this map preserves colimits, which implies that it necessarily is of
the form $X \mapsto f_!(X \otimes D_f)$ for a unique $D_f \in \Sp^T$. 

We recall from~\cite[Theorem~5.6.2.10]{lurieha} that every group in
anima $G$ is the loop group $\Omega BG$ of a pointed connected anima
$$\xymatrix{
  { 1 } \ar[r]^-{s} &
  { BG, } \cr
}$$
and we first apply the general theory above to the unique map
$$\xymatrix{
  { BG } \ar[r]^-{f} &
  { 1. } \cr
}$$
If $G$ is the group in anima underlying a compact Lie group $G$, then
$D_f \simeq S^{\mathfrak{g}}$ is the suspension spectrum of the
one-point compactification of its Lie algebra with the adjoint
$G$-action. So in this situation, the defining fiber sequences takes
the form
$$\xymatrix{
  { f_!(S^{\mathfrak{g}} \otimes X) } \ar[r] &
  { f_*(X) } \ar[r] &
  { f_*^T(X), } \cr
}$$
where $X \in \Sp^{BG}$. This sequence is commonly written as
$$\xymatrix{
  { (S^{\mathfrak{g}} \otimes X)_{hG} } \ar[r] &
  { X^{hG} } \ar[r] &
  { X^{tG}. } \cr
}$$
The Postnikov filtration of $X$ gives rise to the spectral sequences
$$\begin{aligned}
  E_{i,j}^2 & = H^{-i}(BG,\pi_j(s^*(X))) \Rightarrow
\pi_{i+j}(X^{hG}) \cr
E_{i,j}^2 & = \hat{H}^{-i}(BG,\pi_j(s^*(X))) \Rightarrow
\pi_{i+j}(X^{tG}), \cr
\end{aligned}$$
and the edge homomorphism of the former,
$$\xymatrix{
  { \pi_j(X^{hG}) } \ar[r]^-{\theta} &
  { \pi_j(s^*(X)), } \cr
}$$
is induced by the unit map
$$\xymatrix{
  { f_*(X) } \ar[r]^-{\theta} &
  { f_*s_*s^*(X) \simeq s^*(X). } \cr
}$$
The groups in the two $E^2$-terms can be interpreted as the cohomology
and Tate cohomology of $BG$ with coefficients in the local system
$\pi_j(s^*(X))$. If $G$ is finite, then these can be identified with
the group cohomology and Tate cohomology of the group $G$ with 
coefficients in the $G$-module $\pi_j(s^*(X))$.

It follows from~\cite[Theorem~I.4.1]{nikolausscholze} that if
$E \in \CAlg(\Sp^{BG})$ is a commutative algebra in spectra with
$G$-action, then the map $f_*(E) \to f_*^T(E)$ promotes to a map of
commutative algebras in spectra. This implies that the spectral
sequences become spectral sequences of bigraded anticommutative rings
and that the edge homomorphism $\theta$ becomes an map of graded
anticommutative rings; see~\cite{hedenlundrognes}. 
If the $G$-action on $E$ is trivial in the sense that $E \simeq
f^*s^*(E)$, then a choice of trivialization and the unit map combine
to give a section
$$\xymatrix{
  { s^*(E) } \ar[r]^-{\sigma} &
  { f_*f^*s^*(E) \simeq f_*(E) } \cr
}$$
of the edge homomorphism $\theta$. So in this situation, the map
$f_*(E) \to f_*^T(E)$ becomes a map of commutative algebras in
$s^*(E)$-modules in spectra. So the induced map on homotopy
groups becomes a map of anticommutative graded
$\pi_*(s^*(E))$-algebras, and the spectral sequences above become
spectral sequences of anticommutative bigraded $\pi_*(s^*(E))$-algebras.

We are interested in the group $G \simeq U(1)$, and since it is
abelian, the $G$-action on its Lie algebra is trivial. So in this case,
the fiber sequence above becomes
$$\xymatrix{
  { \Sigma E_{hG} } \ar[r] &
  { E^{hG} } \ar[r] &
  { E^{tG}. } \cr
}$$
For instance, in the case of $E \simeq \HH(A/R)$, this gives Connes'
sequence
$$\xymatrix{
  { \Sigma \HC(A/R) } \ar[r] &
  { \HC^{-}(A/R) } \ar[r] &
  { \HP(A/R) } \cr
}$$
relating cyclic homology, negative cyclic homology, and periodic
cyclic homology. Choosing a generator $\bar{v} \in H^2(BU(1),\mathbb{Z})$,
the spectral sequences take the form\footnote{\,Here $\pi_*(-)$
  indicates homotopy groups rather than pushforward along some map $\pi$.}
$$\begin{aligned}
  E^2 & = \pi_*(s^*(E))[\bar{v}] \Rightarrow \pi_*(E^{hU(1)}) \cr
  E^2 & = \pi_*(s^*(E))[\bar{v}^{\pm1}] \Rightarrow \pi_*(E^{tU(1)}). \cr
\end{aligned}$$
We refer to the filtrations of the respective abutments induced by the two
spectral sequences as the ``Nygaard'' filtrations.

We will consider $E \in \CAlg(\Sp^{BU(1)})$ equipped with a
``Bott'' element
$$\xymatrix{
  { \beta \in \pi_2(E^{hU(1)}) } \ar[r] &
  { \pi_2(E^{tU(1)}) } \cr
}$$
and obtain the desired filtered ring from the filtered graded ring
$\pi_*(E^{tU(1)})$ by the following procedure:
\begin{enumerate}
\item[(1)]Invert the Bott element.
\item[(2)]Complete with respect to the Nygaard filtration.
\item[(3)]Extract the filtered subring consisting of homogeneous elements of
  degree $0$.
\end{enumerate}
The image of the induced ring homomorphism
$$\xymatrix{
  { (\pi_*(E^{hU(1)})[\beta^{-1}]^{\wedge})_0 } \ar[r] &
  { (\pi_*(E^{tU(1)})[\beta^{-1}]^{\wedge})_0 } \cr  
}$$
agrees with the $0$th stage of the Nygaard filtration. If
$\pi_*(s^*(E))$ is concentrated in even degrees, then so are
$\pi_*(E^{hU(1)})$ and $\pi_*(E^{tU(1)})$ and all differentials in the
respective spectral sequences are zero. In this case, the map above
is injective.

Suppose first that the $U(1)$-action on $E$ is trivial. If the
homotopy groups $\pi_*(s^*(E))$ are concentrated in even degrees, then
there are isomorphisms
$$\begin{aligned}
  \pi_*(E^{hU(1)}) & \simeq \textstyle{ \pi_*(s^*(E))[[v]]
    \simeq \lim_n \pi_*(s^*(E))[v]/(v^{n+1}) }\cr
  \pi_*(E^{tU(1)}) & \simeq \textstyle{ \pi_*(s^*(E))((v))
    \simeq \pi_*(s^*(E))[[v]][v^{-1}] } \cr
\end{aligned}$$
of anticommutative graded $\pi_*(s^*(E))$-algebras. Here the limit is
calculated in the category of graded $\pi_*(s^*(E))$-algebras and
$v \in \pi_{-2}(E^{hU(1)})$ is a choice of lift of $\bar{v}$, that
is, a ``complex orientation'' of $E$. So with $t^{-1} =
\beta v$, we have
$$\begin{aligned}
  (\pi_*(E^{hU(1)})[\beta^{-1}])_0 & \simeq
  (\pi_*(s^*(E))[\beta^{-1}])_0[[t^{-1}]] \cr 
  (\pi_*(E^{tU(1)})[\beta^{-1}])_0 & \simeq
  (\pi_*(s^*(E))[\beta^{-1}])_0((t^{-1})) \cr
\end{aligned}$$
Hence, writing $R =
(\pi_*(s^*(E))[\beta^{-1}])_0$, we see that the graded
$\pi_*(s^*(E))$-algebra $\pi_*(E^{tU(1)})$ with the Nygaard filtration
gives rise to the $R$-algebra
$$R((t^{-1}))$$
with the $t$-adic filtration. This is the filtered $R$-algebra
that we would obtain from the formal neighborhood of an $R$-valued
point $\infty \in |X|$ of a curve $X \to \Spec(R)$.

Let us now consider the Fargues--Fontaine curve
$$\xymatrix{
  { X = X_F } \ar[r] &
  { \Spec(\mathbb{Q}_p) } \cr
}$$
associated to an algebraically closed completely valued field
$(F,\mathcal{O}_F)$ of characteristic $p > 0$, and let $\infty \in
|X|$ be a closed point corresponding to an untilt
$((C,\mathcal{O}_C),\iota)$ of $(F,\mathcal{O}_F)$. We consider the
commutative algebra in spectra with $U(1)$-action
$$E = \THH(\mathcal{O}_C,\mathbb{Z}_p) \in \CAlg(\Sp^{BU(1)})$$
given by the $p$-adic completion of the topological Hochschild
homology of $\mathcal{O}_C$. In this case, Bhatt--Morrow--Scholze show
in~\cite[Proposition~6.2]{bhattmorrowscholze2} that
$$\begin{aligned}
\pi_*(E^{hU(1)}) & \simeq \TC_*^{-}(\mathcal{O}_C,\mathbb{Z}_p) \simeq
A_{\inf}(\mathcal{O}_F)[u,v]/(uv - \xi) \cr
\pi_*(E^{tU(1)}) & \simeq \TP_*(\mathcal{O}_C,\mathbb{Z}_p) \simeq
A_{\inf}(\mathcal{O}_F)[v^{\pm1}],
\end{aligned}$$
where $u \in \TC_2^{-}(\mathcal{O}_C,\mathbb{Z}_p)$ and
$v \in \TC_{-2}^{-}(\mathcal{O}_C,\mathbb{Z}_p)$, and where $\xi$ is a
generator of the kernel of the edge homomorphism
$$\xymatrix{
  { \TC_0^{-}(\mathcal{O}_C,\mathbb{Z}_p) \simeq
    A_{\inf}(\mathcal{O}_F) } \ar[r]^-{\theta} &
  { \mathcal{O}_C. } \cr
}$$
The element $u$ is called a B\"{o}kstedt element. It is not a Bott
element, but rather a divided Bott element. More precisely, the
element
$$\beta = \varphi^{-1}(\mu)u$$
is a Bott element, but $\varphi^{-1}(\mu) \in A_{\inf}(\mathcal{O}_F)$
is not a unit; see also~\cite[Section~3.2]{hesselholtnikolaus}. So
by inverting $\beta$, we also invert $\varphi^{-1}(\mu)$, which, in
particular, inverts $p$, since
$$\theta(\varphi^{-1}(\mu)) = \zeta_p - 1 \in \mathcal{O}_C$$
is a pseudo-uniformizer. It follows that
$$(\pi_*(E^{tU(1)})[\beta^{-1}]^{\wedge})_0
\simeq (\TP_*(\mathcal{O}_C,\mathbb{Z}_p)[\beta^{-1}]^{\wedge})_0
\simeq B_{\dR}$$
as a filtered $\mathbb{Q}_p$-algebra. It is a discrete valuation field
with residue field $C$. We stress that, after inverting $\beta$, we do not
subsequently $p$-complete, which would leave us with zero. So this
procedure is distinct from Morava $K(1)$-localization.

\section{The affine complement of $\infty \in |X|$}

We will next explain how to obtain the right-hand vertical map
in~(\ref{eq:farguesfontainecurve}) from the ``crystalline Chern
character'' of
Antieau--Mathew--Morrow--Nikolaus~\cite{antieaumathewmorrownikolaus}. The
definition of this map relies on the modern approach to topological
cyclic homology and cyclotomic spectra due to
Nikolaus--Scholze~\cite{nikolausscholze}.

To briefly recall the definition, let $p \colon BU(1) \to BU(1)$ be the
map of pointed anima induced by the $p$-power map of the abelian group
$U(1)$, and let
$$\xymatrix{
  { \Sp_p^{BU(1)} } \ar[r]^-{p_*^T} &
  { \Sp_p^{BU(1)} } \cr
}$$
be the associated Tate cohomology functor. A $p$-complete cyclotomic
spectrum is a pair $(X,\varphi)$ of a $p$-complete
spectrum\footnote{\, The inclusion $\smash{ j_* \colon \Sp_p \to \Sp
  }$ of the full subcategory of $p$-complete spectra admits a left
  adjoint $\smash{ j^* \colon \Sp \to \Sp_p }$, and the unit map
  $\smash{ \eta \colon X \to j_*j^*(X) }$ is $p$-completion. The
  functor $j^*$ admits an essentially unique promotion to
  a symmetric monoidal functor, which, in turn, promotes $j_*$ to a
  lax symmetric monoidal functor. Hence, we obtain an induced
  adjunction of the associated $\infty$-categories of commutative
  algebras.} with $U(1)$-action and a ``Frobenius'' map
$$\xymatrix{
  { X } \ar[r]^-{\varphi} &
  { p_*^T(X) } \cr
}$$
of spectra with $U(1)$-action. Nikolaus and Scholze show that
$p$-complete cyclotomic spectra can be organized into a stable symmetric
monoidal $\infty$-category equipped with a conservative symmetric
monoidal forgetful functor
$$\xymatrix{
  { \CycSp_p } \ar[r] &
  { \Sp_p^{BU(1)} } \cr
}$$
to the symmetric monoidal stable $\infty$-category of $p$-complete
spectra with $U(1)$-action. The tensor unit is given by the pair
$\mathbf{1} \simeq (f^*(\mathbb{S}_p),\varphi)$ of the $p$-complete
sphere spectrum with trivial $U(1)$-action and the composition
$$\xymatrix{
  { f^*(\mathbb{S}_p) } \ar[r] &
  { p_*f^*(\mathbb{S}_p) } \ar[r] &
  { p_*^Tf^*(\mathbb{S}_p) } \cr
}$$
of the adjunct of the equivalence $p^*f^* \simeq (fp)^* \simeq f^*$
and the canonical map. Now, $p$-complete topological
cyclic homology is the symmetric monoidal functor
$$\xymatrix@C=12mm{
  { \CycSp_p } \ar[r]^-{\phantom{,}\TC\phantom{,}} &
  { \Sp_p } \cr
}$$
corepresented by $\mathbf{1} \in \CycSp_p$. If $(X,\varphi)$ is a $p$-complete
cyclotomic spectrum with $X$ bounded below, then $\TC(X,\varphi)$ is
given by the equalizer
$$\xymatrix{
  { \TC(X,\varphi) } \ar[r] &
  { \TC^{-}(X) \simeq f_*(X) } \ar@<.75ex>[r]^-{\varphi}
  \ar@<-.75ex>[r]_-{\phantom{t}\can\phantom{t}} &
  { \TP(X) \simeq f_*^T(X) } \cr
}$$
of the canonical map ``$\can$'' and the cyclotomic Frobenius
``$\varphi$'' defined as follows. Since $fp \simeq f$, we have a
diagram of $p$-complete spectra
$$\begin{xy}
  (3.8,14)*+{ \Sigma f_!p_!(X) \simeq \Sigma f_!(X) }="00";
  (0,14)*+{ \phantom{f_*} }="11";
  (30,14)*+{ f_*(X) }="12";
  (60,14)*+{ f_*^T(X) }="13";
  (0,0)*+{ f_*p_!(X) }="21";
  (30,0)*+{ f_*p_*(X) }="22";
  (60,0)*+{ f_*p_*^T(X) }="23";
  { \ar "12";"00";};
  { \ar "13";"12";};
  { \ar "21";"11";};
  { \ar@{=} "22";"12";};
  { \ar "23";"13";};
  { \ar "22";"21";};
  { \ar "23";"22";};
\end{xy}$$
in which the two rows are fiber sequences. Since $X$ is bounded below,
the Tate orbit lemma~\cite[Lemma~I.2.1]{nikolausscholze} shows
that the left-hand vertical map is an equivalence, and hence, so is the
right-hand vertical map. So in this case, we have the map
$$\xymatrix@C=12mm{
  { X^{hU(1)} \simeq f_*(X) } \ar[r]^-{f_*(\varphi)} &
  { f_*p_*^T(X) \simeq f_*^T(X) \simeq X^{tU(1)}, } \cr
}$$
which, by abuse of notation, we also denote by $\varphi \colon
\TC^{-}(X) \to \TP(X)$.

We recall the most economical way of associating to a commutative
algebra in $p$-complete spectra $R$ a commutative algebra
in $p$-complete cyclotomic spectra
$$(\THH(R,\mathbb{Z}_p),\varphi)$$
following~\cite[Section~IV.2]{nikolausscholze}.

In general, let $G$ be a group in anima, and let $s \colon 1 \to BG$
be the corresponding connected pointed anima. If $\mathcal{C}$ is any
presentable $\infty$-category, then we have
$$\xymatrix{
  { \mathcal{C} } \ar@<.7ex>[r]^-{s_!} &
  { \mathcal{C}^{BG}, } \ar@<.7ex>[l]^-{s^*} \cr
}$$
where the right adjoint $s^*$ takes an object of $\mathcal{C}$
with $G$-action to the underlying object of $\mathcal{C}$, and where the
left adjoint $s_!$ takes an object of $\mathcal{C}$ to the
free object of $\mathcal{C}$ with $G$-action. We have a cartesian
square of pointed anima
$$\xymatrix{
  { G } \ar[r]^-{s'} \ar[d]^-{s'} &
  { 1 } \ar[d]^-{s} \cr
  { 1 } \ar[r]^-{s} &
  { BG } \cr
}$$
and the base-change formula applied to this square shows that
$$s^*s_! \simeq s'_!s'{}^*$$
as endofunctors of $\mathcal{C}$. Informally, if $X$ is an object of
$\mathcal{C}$, then $s^*s_!(X)$ is the object of $\mathcal{C}$
underlying the free object of $\mathcal{C}$ with $G$-action 
associated with $X$, whereas $s'_!s'{}^*(X)$ is the colimit in
$\mathcal{C}$ of the constant diagram with value $X$ indexed by $G$.

We apply this with $G \simeq U(1)$ and $\mathcal{C} \simeq
\CAlg(\Sp_p)$. Given $R \in \CAlg(\Sp_p)$, we define its $p$-complete
topological Hochschild homology to be the commutative algebra in
$p$-complete spectra with $U(1)$-action given by
$$\THH(R,\mathbb{Z}_p) \simeq s_!(R) \in \CAlg(\Sp_p)^{BU(1)}$$
and the base-change formula shows that its underlying commutative
algebra in $p$-complete spectra is given by
$$s^*\THH(R,\mathbb{Z}_p) \simeq s^*s_!(R) \simeq s'_!s'{}^*(R) \simeq
R^{\hspace{.5pt}\otimes\hspace{.5pt} U(1)},$$
where the right-hand term is suggestive notation for the colimit in
$\CAlg(\Sp_p)$ of the constant diagram with value $R$ indexed by
the anima underlying $U(1)$.

We next recall the definition of the Frobenius map
$$\xymatrix{
  { \THH(R,\mathbb{Z}_p) } \ar[r]^-{\varphi} &
  { p_*^T(\THH(R,\mathbb{Z}_p)), } \cr
}$$
where we abuse notation and also write $p_*^T$ for the
endofunctor of
$$\CAlg(\Sp_p)^{BU(1)} \simeq \CAlg(\Sp_p^{BU(1)})$$
induced by the lax symmetric monoidal endofunctor $p_*^T$ of
$\Sp_p^{BU(1)}$. We consider the following diagram of pointed anima 
$$\xymatrix@C=10mm{
  { 1 } \ar[dr]^-{t} \ar@/_.7pc/[ddr]_-{\id} \ar@/^.7pc/[drr]^-{s} &
  {} &
  {} \cr
  {} &
  { BC_p } \ar[r]^-{i} \ar[d]^-{g} &
  { BU(1) } \ar[d]^-{p} \cr
  {} &
  { 1 } \ar[r]^-{s} &
  { BU(1) } \cr
}$$
with the square cartesian. As part of the structure of any symmetric
monoidal $\infty$-category $\mathcal{C}$, there is a symmetric
monoidal functor
$$\xymatrix{
  { \mathcal{C} } \ar[r]^-{t_!^{\otimes}} &
  { \mathcal{C}^{BC_p} } \cr
}$$
which, informally, takes $X$ to $X^{\otimes p}$ together with an
equivalence between the induced functor of commutative algebras 
and
$$\xymatrix{
  { \CAlg(\mathcal{C}) } \ar[r]^-{t_!} &
  { \CAlg(\mathcal{C})^{BC_p} \simeq \CAlg(\mathcal{C}^{BC_p}). } \cr
}$$
We now recall from~\cite[Proposition~III.3.1]{nikolausscholze} that there
is a unique lax symmetric monoidal transformation called the Tate
diagonal\footnote{\,The generalized Segal conjecture for $C_p$ proved
  in~\cite[Theorem~III.1.7]{nikolausscholze} states that if $X$ is
  bounded below, then the Tate diagonal is an equivalence.} 
$$\xymatrix{
  { X } \ar[r]^-{\delta} &
  { g_*^Tt_!^{\otimes}(X). } \cr
}$$
So for $R \in \CAlg(\Sp_p)$, we get the composite map
$$\xymatrix{
  { R } \ar[r]^-{\delta} &
  { g_*^Tt_!(R) } \ar[r]^-{\eta} &
  { g_*^Ti^*i_!t_!(R) \simeq g_*^Ti^*s_!(R) \simeq s^*p_*^Ts_!(R) } \cr
}$$
of commutative algebras in $p$-complete spectra, whose adjunct
$$\xymatrix{
  { s_!(R) } \ar[r]^-{\varphi} &
  { p_*^Ts_!(R) } \cr
}$$
is the desired map of commutative algebras in $p$-complete spectra
with $U(1)$-action.

The crystalline Chern character of~\cite{antieaumathewmorrownikolaus}
is obtained as follows. Given a $p$-complete cyclotomic spectrum $X =
(X,\varphi)$, we may form the tensor product
$$X \otimes \THH(\mathbb{F}_p,\mathbb{Z}_p)$$
in the symmetric monoidal $\infty$-category $\CycSp_p$ of $p$-complete
cyclotomic spectra and consider the canonical map
$$\xymatrix{
  { \TC(X \otimes \THH(\mathbb{F}_p,\mathbb{Z}_p)) } \ar[r] &
  { \TC^-(X \otimes \THH(\mathbb{F}_p,\mathbb{Z}_p)) } \cr
}$$
from its topological cyclic homology to its negative topological
cyclic homology. We consider this map for $X \simeq
\THH(\mathcal{O}_C,\mathbb{Z}_p)$, take homotopy groups, invert the
Bott element, Nygaard complete the target, and extract the subrings of
homogeneous elements of degree
$0$. By~\cite[Theorem~2.12]{antieaumathewmorrownikolaus}, the
resulting map takes the form
$$\xymatrix{
  { (\TC_*(\mathcal{O}_C/p,\mathbb{Z}_p)[\beta^{-1}])_0 } \ar[r] &
  { (\TP_*(\mathcal{O}_C,\mathbb{Z}_p)[\beta^{-1}]^{\wedge})_0, } \cr
}$$
and this map is the desired crystalline Chern
character. Finally, the identification of this map with the map $B_e
\to B_{\dR}$ in~(\ref{eq:farguesfontainecurve}) is a consequence
of~\cite[Corollary~8.23]{bhattmorrowscholze2}
and~\cite[Proposition~6.2]{bhattmorrowscholze2}. We conclude that the
diagram
\begin{equation}\label{eq:tc}
  \xymatrix{
    {} &
    { (\TC_*(\mathcal{O}_C/p,\mathbb{Z}_p)[\beta^{-1}])_0 } \ar[d] \cr
    { \Fil_n(\TP_*(\mathcal{O}_C,\mathbb{Z}_p)[\beta^{-1}]^{\wedge})_0
    }  \ar[r] &
    { (\TP_*(\mathcal{O}_C,\mathbb{Z}_p)[\beta^{-1}]^{\wedge})_0 } \cr
  }
\end{equation}
obtained from the various flavors of topological cyclic homology and
maps between them is canonically identified with the diagram of period
rings~(\ref{eq:farguesfontainecurve}).

\section{Dreams}

We would like to similarly obtain the
diagram~(\ref{eq:twistorprojectiveline}) for the twistor projective
line from a diagram analogous to~(\ref{eq:tc}). Given a closed point
$\infty \in |X|$ with residue field $k(\infty) \simeq \mathbb{C}$, one
might hope to obtain the bottom line
in~(\ref{eq:twistorprojectiveline}) from the periodic cyclic homology
$\HP_*(\mathbb{C}/\hspace{.7pt}\mathbb{R})$. This will not work,
however, due to the lack of a Bott element $\beta \in
\HC_2^{-}(\mathbb{C}/\hspace{.7pt}\mathbb{R})$. The liquid theory of
Clausen--Scholze~\cite{clausenscholzeanalytic} points to an
alternative. Let $0 < r < 1$ be a real number. The Harbater ring is  
the subring
$$\mathbb{Z}((T))_{>r} \subset \mathbb{Z}((T))$$
consisting of the Laurent series that for some $r' > r$ converge
absolutely on a disc of radius $r'$ centered at the origin. Harbater
shows in~\cite{harbater} that this ring is a principal ideal domain; see
also~\cite[Theorem~7.1]{clausenscholzeanalytic}. In particular, if
$x \in \mathbb{C}$ and $|x| \leq r$, then the kernel of 
the continuous ring homomorphism
$$\xymatrix{
  { \mathbb{Z}((T))_{>r} } \ar[r] &
  { \mathbb{C} } \cr
}$$
that to $T$ assigns $x$ is a principal ideal, and hence, the
Hochschild homology
$$\HH_*(\mathbb{C}/\mathbb{Z}((T))_{>r}) = \mathbb{C}\langle \bar{u}
\rangle = \mathbb{C}[\bar{u}]$$
is a (divided) polynomial algebra on a generator $\bar{u}$ of degree
$2$. So with
$$E = \HH(\mathbb{C}/\mathbb{Z}((T))_{>r}) \in \CAlg(\Sp^{BU(1)}),$$
we find that
$$\begin{aligned}
\pi_*(E^{hU(1)}) & \simeq \HC_*^{-}(\mathbb{C}/\mathbb{Z}((T))_{>r})
\simeq \mathbb{C}[[\xi]][u,v]/(uv-\xi) \cr
\pi_*(E^{tU(1)}) & \simeq \HC_*^{-}(\mathbb{C}/\mathbb{Z}((T))_{>r})
\simeq \mathbb{C}[[\xi]][v^{\pm1}] \cr
\end{aligned}$$
with $u \in \HC_2^{-}(\mathbb{C}/\mathbb{Z}((T))_{>r})$ and $v \in
\HC_{-2}^{-}(\mathbb{C}/\mathbb{Z}((T))_{>r})$ and with $\xi$ is a
generator of the kernel of the edge homomorphism
$$\xymatrix{
  { \HC_0^{-}(\mathbb{C}/\mathbb{Z}((T))_{>r}) \simeq
    \mathbb{C}[[\xi]] } \ar[r]^-{\theta} &
  { \mathbb{C}, } \cr
}$$
analogously to $E \simeq \THH(\mathcal{O}_C,\mathbb{Z}_p)$. Hence, for a
``Bott element'' of the form $\beta = fu$ with
$f \in \mathbb{C}[[\xi]]^{\times}$ a unit, we obtain the bottom row
in~(\ref{eq:twistorprojectiveline}) with $t^{-1} = \beta v$.

\section*{Acknowledgements}

The author would like to express his gratitude to an anonymous referee for a
number of helpful comments.

\providecommand{\bysame}{\leavevmode\hbox to3em{\hrulefill}\thinspace}
\providecommand{\MR}{\relax\ifhmode\unskip\space\fi MR }
\providecommand{\MRhref}[2]{%
  \href{http://www.ams.org/mathscinet-getitem?mr=#1}{#2}
}
\providecommand{\href}[2]{#2}

\end{document}